\documentclass[%
preprint,
showpacs,
 amsmath,amssymb,
 aps,
pra,
]{revtex4-1}

\usepackage{graphicx}
\usepackage{dcolumn}
\usepackage{bm}

\begin{document}


\title{Ro-Vibrational Hamiltonian of Three Body Systems Near Collinear Configurations}

\author{\"{U}nver \c{C}ift\c{c}i}
\affiliation{Department of Mathematics, Nam\i k Kemal University, 59030 Tekirda\u{g}, Turkey}

\date{\today}

\begin{abstract}
Recent developments on three body systems have revealed that dynamics of trajectories passing through collinear configurations can be easily adopted. We analyse the reduction procedure in order to detect the points where collinear configurations are deviating. Then we show that the value of the reduced Hamiltonian can be computed at these points.

\end{abstract}

\maketitle


\section{\label{sec:level1}Introduction}

Reduction of three body systems dates back to Lagrange. His reduction can be stated as using invariants in general terms. Choosing those invariants is the key difference among reduction methods.    

Lagrange's original reduction uses mutual distances as coordinates, so the invariance under  rotations is ensured \cite{AlbouyChenciner98}. There are ten coordinates and two Casimirs in this reduction method. But the reduced space is eight dimensional and in this method obtaining physical or effective coordinates seems to be impossible. Also partitioning the dynamics into internal and rotational is not evident. Nevertheless this method has several advantages, for instance it respects collinear configurations. This method is also useful for reducing boost symmetries for the gravitational three body systems.

Jacobi's elimination of nodes is another method of reduction which enables to use internal variables so that the dynamics can be partitioned into internal and rotational parts. It relies on the fact that when describing dynamics in body frames, internal dynamics is partitioned as much as possible.  A geometric picture of Jacobi's elimination of nodes has recently been made evident \cite{Iwai87a, LittlejohnReinsch97, IwaiYamaoka08}. So far motions passing through collinear configurations hasn't been considered in this picture.

We review below some of the other methods of reduction:

A splitting method is used Moeckel \cite{Moeckel04} to study dynamics near a collinear relative equilibrium, equilibrium of the reduced Hamiltonian. Choosing a suitable basis, it is possible to linearize the equations motion around the collinear equilibrium. In this method,  the Hamiltonian is not obtained and the coordinates are not used.

Hsiang's regularization is another way of treating possible collinear configurations along a motion of three bodies \cite{Montgomery96}. This method  uses oriented triangles. It is useful for the Newtonian problem. But obtaining the reduced Hamiltonian seems to be impossible. 

In physics literature so-called isomorphic Hamiltonian is introduced to unify the Hamiltonians of collinear and noncollinear configurations \cite{Watson70}. A non-physical object is attached and that introduces an extra angle parameter. But it is not considered that collinear three body systems have two effective degrees of freedom in this method.   

We make use of the fact that that the action of rotations on the translation reduced state space is free. This makes identification of a well-defined body frame possible. We use the space fixed angular momentum to define the body frame. The body frame is shown to limit to the frame at the collinear configuration.

As the reduced Hamiltonian can be defined at noncollinear configurations, and as it has a limit at a collinear configuration, we just need to find its value at the collinear configuration. We first compute the value of the Hamiltonian by using the Legendre transform. Then we show that this value is the limit of the noncollinear Hamiltonian at the collinear configuration.

Reduced Hamiltonian obtained by Jacobi's elimination of nodes which is equivalent to ro-vibrational Hamiltonian in Physics literature is important in the sense that it has the internal and rotational partitioning as much as possible. Its drawback was that it it couldn't be defined at collinear configurations where three bodies lie on a line. We show that the ro-vibrational Hamiltonian has a well-defined value at these points which is its limit when the motion tends to a collinear configuration. The value is computed in a way that a body frame is defined by the use of angular momentum vector. This ensures that the body frame at the collinear configuration is the limit of body frames when the motion tends to the collinear configuration.

Ro-vibrational Hamiltonian is widely used in molecular physics and related areas \cite{LittlejohnReinsch97}. Although it is very effective for the study of dynamics of three and more bodies, it's been not possible to define it for motions passing through collinear configurations. Defining this Hamiltonian at collinear shapes is crucial in both theoretical and practical point of view. 

\subsection*{Three body Reduction}

Consider a  system of three bodies with masses $m_1,m_2,m_3$ and position vectors $\mathbf{x}_{1},\mathbf{x}_{2},\mathbf{x}_{3} \in \mathbb{R}^3$,
respectively, with potential energy which is invariant under translations and rotations, and without external forces acting on the three bodies. The symmetry of overall translations can be reduced, for instance, by introducing
mass-weighted Jacobi vectors which are defined according to  
\begin{eqnarray*}
\mathbf{s}_1&=&\sqrt{\mu_1}(\mathbf{x}_{1}-\mathbf{x}_{3}),\\
\mathbf{s}_2&=&\sqrt{\mu_2} 
(\mathbf{x}_{2}-\frac{m_{1}\mathbf{x}_{1}+m_{3}\mathbf{x}_{3}}{m_{1}+m_{3}}),
\end{eqnarray*}%
where
\begin{equation*}
 \mu_1=\frac{m_{1}m_{3}}{m_{1}+m_{3}}, \ \ \ \mu_2=\frac{m_{2}(m_{1}+m_{3})}{m_{1}+m_{2}+m_{3}}
\end{equation*}
are  reduced masses. 

Action of rotations on the state space is free for generic configurations. To see this, let the angular momentum 
$\mathbf{L}$ defined by
\begin{equation}
 \mathbf{L}=\mathbf{s}_1 \times \dot{\mathbf{s}}_1 + \mathbf{s}_2 \times \dot{\mathbf{s}}_2
\end{equation}
be nonzero. Then clearly $\dot{\mathbf{s}}_1$ or $\dot{\mathbf{s}}_2$ is not parallel to $\mathbf{s}_1$  and $\mathbf{s}_2$. So the Lagrangian and the corresponding Hamiltonian can be reduced on a well-defined quotient space. The reduction is usually done by using invariant parameters which are physical quantities of the system. Here we review the most used one in Physics which is called ro-vibrational Hamiltoinan, with a special care at collinear configurations.   

We chose a space frame $\mathbf{e}_1, \mathbf{e}_2, \mathbf{e}_3$ and a body frame $\mathbf{u}_1, \mathbf{u}_2, \mathbf{u}_3$. Expressing dynamics in terms of the body frame gives the reduction. At each instant $t$ there is a rotation $\mathsf{R}(t)$ which rotates the space frame to the body frame, so if we choose the Euler angles by 
\begin{widetext}
\begin{eqnarray}
\mathbf{u}_1 &=& \sin \beta \cos \alpha \, \mathbf{e}_1 + \sin \beta \sin \alpha \, \mathbf{e}_2 + \cos \beta \, \mathbf{e}_3, \nonumber \\
\mathbf{u}_2 &=& (- \cos \beta \cos \alpha \sin \gamma - \sin \alpha \cos \gamma) \, \mathbf{e}_1 + (- \cos \beta \sin \alpha \sin \gamma + \cos \alpha \cos \gamma) \, \mathbf{e}_1 +  \sin \beta \sin \gamma \,\mathbf{e}_3,  \\
\mathbf{u}_3 &=& (\cos \beta \cos \alpha \cos \gamma - \sin \alpha \sin \gamma) \, \mathbf{e}_1 + (\cos \beta \sin \alpha \cos \gamma + \cos \alpha \sin \gamma) \, \mathbf{e}_2 - \sin \beta \sin \gamma \,\mathbf{e}_3, \nonumber \label{euler_angles}
\end{eqnarray}
\end{widetext}
where $0 \le \alpha \le 2 \pi, 0 \le \beta \le \pi, 0 \le \gamma \le 2 \pi$ we have a matrix representation of $\mathsf{R}$. 

Let us denote the representation of the Jacobi vectors in the space frame by $\mathbf{s}_1, \mathbf{s}_2$ as above and let the representation of there vectors in body frame be denoted by $\mathbf{r}_1, \mathbf{r}_2$. So they are related by
 \begin{equation}
 \mathbf{s}_i=\mathsf{R} \, \mathbf{r}_i, \hspace{3mm} i=1,2\,.
\end{equation}
For noncollinear shapes this equation fully specifies  $\mathsf{R}$ but for the collinear ones this not the case. We will show that for a motion passing through a collinear configuration the body frame can be identified as the limit of the body frames of the noncollinear shapes.  Let the Jacobi coordinates $(r_1,r_2, \phi)$ of the Jacobi vectors in body frame be defined by
\begin{eqnarray*}
\mathbf{r}_1 &=&r_1 \, \mathbf{u}_{1}, \\
\mathbf{r}_2 &=&r_2\, \cos \phi \, \mathbf{u}_{1}+r_2\sin \phi \,  \mathbf{u}_{2},
\end{eqnarray*}%
where $0 \le \phi \le \pi$. We assume that the collinear configuration of a motion corresponds to $t=0$ and $\phi(0) =0$. 

We note here that if the motion is fully collinear, then two Euler angles $\alpha$ and $\beta$ (which can be seen by the definition of the Euler angles) and two internal coordinates $r_1$ and $r_2$ can be used to study the reduced dynamics \cite{IwaiYamaoka05a}.

Now we  want to obtain the value of the reduced Hamiltonian at the collinear configuration. For this purpose we are carefully adopting the classical procedure.
The kinetic energy is given by
\begin{equation}
K=\frac{1}{2} \sum^{3}_{i=1} m_i \, \mathbf{\dot{x}}^{2}_i\,=\frac{1}{2} \sum^{2}_{i=1} \dot{\mathbf{s}}_i^{2}. \,
\end{equation}
Defining body velocities according to
\begin{equation}
 \mathbf{v}_i=\mathsf{R}^T \, \dot{\mathbf{s}}_i, \hspace{3mm} i=1,2,
\end{equation}
and using the shape coordinates and their time derivatives one can rewrite the body velocities as
\begin{eqnarray*}
\mathbf{v}_i&=&\mathsf{R}^T \, (\dot{\mathsf{R}} \, \mathbf{r}_i  
+\sum^{3}_{\mu=1} \mathsf{R} \, \frac{\partial \mathbf{r}_i}{\partial q^{\mu}}\, \dot{q}^\mu) \\
&=&\mathsf{R}^T \, \dot{\mathsf{R}} \, \mathbf{r}_i +\sum^{3}_{\mu=1} \frac{\partial \mathbf{r}_i}{\partial q^{\mu}}\, \dot{q}^\mu,
\end{eqnarray*}
where $q^1=r_1, q^2=r_2, q^3=\phi$.
The body angular velocity $\boldsymbol{\omega}=(\omega_1,\omega_2,\omega_3)$ is the vector
in $\mathbf{R}^3$ corresponding to the skew-symmetric matrix $\mathbf{\Omega}=\mathsf{R}^T \, \dot{\mathsf{R}}$
given by the identification
\begin{equation}
\left(
\begin{array}{ccc}
0 & -\omega_3 & \omega_2 \\ 
\omega_3 & 0 & -\omega_1 \\ 
-\omega_2 & \omega_1 & 0
\end{array}
\right) \leftrightarrow
\left(
\begin{array}{ccc}
\omega_1 \\ 
\omega_2  \\ 
\omega_3 
\end{array}
\right)\,. 
\end{equation}
Then one has
\begin{equation}
\mathbf{v}_i=\boldsymbol{\omega} \times \mathbf{r}_i+\sum^{3}_{\mu=1} \frac{\partial \mathbf{r}_i}{\partial q^\mu} \, \dot{q}^\mu. \label{velocity}
\end{equation}

In terms of Euler angles one can compute the components of $\boldsymbol{\omega}$ as%
\begin{eqnarray*}
\omega _{1}&=&\dot{\alpha} \, \cos \beta + \dot{\gamma}  ,\\
\omega _{2}&=& \dot{\beta} \, \cos \gamma  +\dot{\alpha} \, \sin \beta \sin \gamma ,\\
\omega _{3}&=& \dot{\beta} \, \sin \gamma  -\dot{\alpha} \, \sin \beta \cos \gamma\, . 
\end{eqnarray*}%

Then the kinetic energy  becomes
\begin{equation}
K=\frac{1}{2} \, \boldsymbol{\omega}^T \, \mathbb{I} \, \boldsymbol{\omega}+\sum^{3}_{\mu=1} (\boldsymbol{\omega}^T \,  \mathbf{a_\mu}) \,\dot{q}^\mu
+\frac{1}{2} \, \sum^{3}_{\mu,\nu=1}h_{\mu \nu}\, \dot{q}^\mu \, \dot{q}^\nu,   
\end{equation}
where 
\begin{equation}
h_{\mu \nu}=\sum^{2}_{i=1} \frac{\partial \mathbf{r}_i}{\partial q^{\mu}}^T \, \frac{\partial \mathbf{r}_i}{\partial q^\nu},
\end{equation}
$\mathbb{I}$ is the moment of inertia tensor given by 
\begin{equation}
 \mathbb{I} \, \mathbf{u}=\mathbf{r}_1 \times (\mathbf{u} \times \mathbf{r}_1)+\mathbf{r}_2 \times (\mathbf{u} \times \mathbf{r}_2),
\end{equation}
for $\mathbf{u}\in \mathbb{R}^3$, 
and 
\begin{equation}
 \mathbf{a_{\mu}}=\mathbf{r}_1 \times \frac{\partial \mathbf{r}_1}{\partial q^\mu}+
\mathbf{r}_2 \times \frac{\partial \mathbf{r}_2}{\partial q^\mu}, \quad \mu=1,2.
\end{equation}

One can specifically obtain
\cite{LittlejohnReinsch97}
\begin{equation}
\mathbb{I}=\left[
\begin{array}{ccc} 
r_{2}^{2} \sin^2 \phi & -r_{2}^{2} \sin \phi \cos \phi & 0 \\ 
-r_{2}^{2} \sin \phi \cos \phi & r_{1}^{2}+r_{2}^{2} \cos^2  \phi & 0 \\ 
0 & 0 & r_{1}^{2}+r_{2}^{2} 
\end{array}
\right]\,, \label{eq:inertia_tensor_Jacobi}
\end{equation}
\begin{equation}
 \left[ h_{\mu\nu}\right] =\left[ 
\begin{array}{ccc}
1 & 0 & 0 \\ 
0 & 1 & 0 \\ 
0 & 0 & r_{2}^{2}%
\end{array}%
\right]\,,
\end{equation}
and
\begin{equation}
 \mathbf{a}_{r_1 }=\mathbf{a}_{r_2}=(0,0,0),\ \ \ \mathbf{a}_{\phi }=(0,0,r_2^2), \label{eq:mech_connection_Jacobi}
\end{equation}
respectively.

The body angular momentum $\mathbf{J}= (J_1,J_2,J_3)$ is given by
\begin{equation}
 \mathbf{J}=\mathsf{R}^T \, \mathbf{L}= \mathbf{r}_1 \times \mathbf{v}_1+\mathbf{r}_2 \times \mathbf{v}_2.
\end{equation}
Then by Equation \eqref{velocity} and one has
\begin{equation}
 \mathbf{J}= \mathbb{I} \, \boldsymbol{\omega}+\sum^{3}_{\mu=1} \mathbf{a}_{\mu} \,\dot{q}^\mu \label{angular}
\end{equation}
or in coordinates
\begin{eqnarray}\label{angular_momentum_components}
J_1&=&r_2^2 \, \omega_1 \, \sin^2 \phi - r_2^2 \, \omega_2 \, \sin \phi \, \cos \phi, \nonumber \\
J_2 &=&  -r_2^2 \, \omega_1 \, \sin \phi \, \cos \phi - (r_1^2+ r_2^2  \, \cos^2 \phi) \, \omega_2,  \\ 
J_3&=&(r_1^2+r_2^2) w_2+r_2^2 \, \dot{q}^3 \,. \nonumber
\end{eqnarray}

At the collinear configuration the kinetic energy takes the form
\begin{widetext}
\begin{equation}
K(0)=\frac{1}{2} (r_1(0)^2+r_2(0)^2)(\omega_2(0)^2+\omega_3(0)^2)+r_2^2 \omega_3(0) \dot{q}^3(0)
+\frac{1}{2} (\dot{q}^1(0))^2+\frac{1}{2} (\dot{q}^2(0))^2+\frac{1}{2}r_2^2  (\dot{q}^3(0))^2.
\end{equation} 
\end{widetext}
We set $\omega_2(0)=0$ for simplicity. This can be achieved, for instance, by assuming that $\alpha(0)=\beta(0)=0$ and $\gamma(0)=\pi/2$ in the convention we use. In this case $\mathbf{u}_1(0)=\mathbf{e}_3, \mathbf{u}_2(0)=-\mathbf{e}_1, \mathbf{u}_3(0)=\mathbf{e}_2 \, $.

Clearly by Legendre transformation, 
\begin{eqnarray}
J_1(0)&=&0, \nonumber \\
J_2(0)&=&0, \\
J_3(0)&=&\frac{\partial K}{\partial \omega_3}\Bigr|_{\substack{t=0}}
=(r_1(0)^2+r_2(0)^2) \, \omega_3(0) + r_2(0)^2 \, \dot{q}^3(0), \nonumber
\end{eqnarray}
and
\begin{eqnarray}
p_1(0)&=&\frac{\partial K}{\partial \dot{q}^1}\Bigr|_{\substack{t=0}}
=\dot{q}^1(0),  \nonumber \\ 
p_2(0)&=&\frac{\partial K}{\partial \dot{q}^2}\Bigr|_{\substack{t=0}}=\dot{q}^2(0), \\
p_3(0)&=&\frac{\partial K}{\partial \dot{q}^3}\Bigr|_{\substack{t=0}}=r_2(0)^2 \, (w_3(0) +\dot{q}^3(0)) \,. \nonumber
\end{eqnarray}
So by the the definition of the Hamiltonian 
\begin{equation}
H(0)=\omega(0)^T \mathbf{J}(0) + \sum_{\mu=1}^3 p_{\mu}(0) \, \dot{q}^{\mu}(0) - L(0) 
\end{equation} 
we have
\begin{widetext} 
\begin{equation}
H(0)=\omega_3(0) \, J_3(0) + p_1(0)^2 + p_2(0)^2+ p_3(0) \, \dot{q}^3(0) - K(0) +V(0) \,. 
\end{equation} 
\end{widetext}
Here one can see that
\begin{equation}
\omega_3(0)=\frac{J_3(0)}{r_2(0)^2+r_2(0)^2} 
\end{equation} 
and 
\begin{equation}
\dot{q}^3(0)=\frac{r_1(0)^2+r_2(0)^2}{r_2(0)^2 \, r_2(0)^2} \, p_3(0) - \frac{1}{r_1(0)^2} \, J_3(0) \,. 
\end{equation} 
Substituting and manipulating we have
\begin{widetext}
\begin{equation}
H(0)= \frac{1}{r_1(0)^2+r_2(0)^2} \, J_3(0)^2  + p_1(0)^2+
p_2(0)^2 +   \frac{r_1(0)^2+r_2(0)^2}{r_1(0)^2 \, r_2(0)^2} \, (p_3(0)-\frac{r_2(0)^2}{r_1(0)^2+r_2(0)^2} \, J_3(0))^2 +V(0) \,. 
\end{equation} 
\end{widetext}

Now we show that the value $H(0)$ is the limit of the reduced Hamiltonian of noncollinear configurations. 
For noncollinear configurations, 
using that the horizontal metric 
\begin{equation}
g_{\mu \nu}=h_{\mu \nu}-\mathbf{A}_{\mu}^T \, \mathbb{I} \, \mathbf{A_\nu},
\end{equation}
where 
\begin{equation}
 \mathbf{A_{\mu}}=\mathbb{I}^{-1} \, \mathbf{a}_{\mu}
\end{equation}
allows one to write  the kinetic energy in the  compact form
\begin{widetext}
\begin{equation}
 K=\frac{1}{2} \, (\boldsymbol{\omega}+\sum^{3}_{\mu=1} \mathbf{A_\mu} \, \dot{q}^\mu)^T \, \mathbb{I} \, (\boldsymbol{\omega}+ \sum^{3}_{\nu=1}
\mathbf{A}_\nu \, \dot{q}^\nu)+\frac{1}{2} \sum^{3}_{\mu,\nu=1} \,g_{\mu \nu} \, \dot{q}^\mu \, \dot{q}^\nu.
\end{equation}
\end{widetext}

The conjugate momenta are given by
\begin{equation}
 \mathbf{J}=\frac{\partial K}{\partial \boldsymbol{\omega}}=\mathbb{I} \, (\mathbf{\xi} +\sum^{3}_{\mu=1} \mathbf{A}_\mu \, \dot{q}^\mu),
\end{equation}
and
\begin{equation}
 p_\mu =\frac{\partial K}{\partial \dot{q}^\mu}=\sum^{3}_{\nu=1}g_{\mu \nu} \, \dot{q}^\nu + \mathbf{J}^T \, \mathbf{A}_\mu. \label{three_moment}
\end{equation}
Thus the Hamiltonian takes the form
\begin{widetext}
\begin{equation}
H=\frac{1}{2}\,\mathbf{J}^T \, \mathbb{I}^{-1} \, \mathbf{J}+
\frac{1}{2} \sum^{3}_{\mu,\nu=1} g^{\mu \nu} (p_{\mu}-\mathbf{J}^T \, \mathbf{A}_\mu)(p_\nu-\mathbf{J}^T \, \mathbf{A}_\nu)+V(q_1,q_2,q_3)  \, . \label{Ham}
\end{equation}
\end{widetext}
Putting the results above together the Hamiltonian in terms of Jacobi coordinates becomes
\begin{widetext}
\begin{eqnarray*}
 H(r_1,r_2,\phi,p_1,p_2,p_3,\mathbf{J})&=&\frac{1}{2} \, \{ \frac{r_{1}^{2}+r_{2}^{2} \cos^2  \phi}{r_{1}^{2} \, r_{2}^{2} \, \sin^2  \phi} \, J_{1}^{2}+
\frac{2 \, \cos \phi}{r_{1}^{2} \, \sin \phi} \, J_{1} \, J_{2}+\frac{1}{r_{1}^{2}}\,J_{2}^{2}+\frac{1}{r_{1}^{2}+r_{2}^{2}}\,J_{3}^{2} \\
&+&p_{1}^{2} +p_{2}^{2}+\frac{r_1^2+r_2^2}{r_1^2 \, r_2^2} \, (p_{3}-\frac{r_{2}^{2}}{r_{1}^{2}+r_{2}^{2}} \,J_{3})^{2} \}+V(r_1,r_2,\phi)\,.
\end{eqnarray*}
\end{widetext}

We show that the noncollinear Hamiltonian limits to the value we computed as the trajectory tends to the collinear configuration. We consider the expression
\begin{equation}
\frac{1}{\sin \phi}  \, J_1 \label{singular}
\end{equation}
which occurs in the first two summands of the final Hamiltonian. By \eqref{angular_momentum_components} it can be written as
\begin{equation}
\frac{1}{\sin \phi}  \,  (r_2^2 \, \omega_1 \, \sin^2 \phi - r_2^2 \, \omega_2 \, \sin \phi \, \cos \phi) \, .\label{singular_two}
\end{equation}
As $t$ tends to zero, this expression tends to  
\begin{equation}
 r_2(0)^2 \, \omega_1(0) \, \sin \phi(0) - r_2(0)^2 \, \omega_2(0) \, \cos \phi(0) \, \label{}
\end{equation}
or 
\begin{equation}
 - r_2(0)^2 \, \omega_2(0) \,  ,\label{}
\end{equation}
since we assume $\phi(0)=0$. We specifically had set $\omega_2(0)=0$ by making suitable assumptions, hence the expression \eqref{singular} tends to zero, as $t$ tends to zero. So, the Hamiltonian limits to the value we computed for the collinear configuration. 

The angular momentum $\mathbf{L}$ is constant along the motion. Consider the molecular line $l$ at the collinear configuration. Apparently $\mathbf{L}$ is orthogonal to the line $l$. We can specify the space frame $\mathbf{e}_1,\mathbf{e}_2,\mathbf{e}_3$ by aligning $\mathbf{e}_2$ with $\mathbf{L}$ and $\mathbf{e}_3$ with $l$. If the body frame is defined as in \eqref{euler_angles}, then assuming $\alpha(0)=\beta(0)=0$ and $\gamma(0)=\pi/2$ results in $\mathbf{u}_1(0)=\mathbf{e}_3, \mathbf{u}_2(0)=-\mathbf{e}_1$ and $\mathbf{u}_3(0)=\mathbf{e}_2$. So that we have identified the body frame uniquely at the collinear configuration. This means that the body fixed frame $\mathbf{u}_1(t),\mathbf{u}_2(t),\mathbf{u}_3(t)$ tends to $\mathbf{u}_1(0),\mathbf{u}_2(0),\mathbf{u}_3(0)$ as $t$ tends to zero.  

\section*{Conclusions}

Dynamics of three body systems used to be studied separately for collinear and noncollinear configurations because the usual procedure makes heavily use of the inverse of moment of inertia tensor which is singular at collinear configurations. We give a simple and algorithmic derivation of the value of the reduced Hamiltonian at a collinear configuration in order to extend the classical ro-vibrational Hamiltonian. We believe that it could be insightful in studying molecular systems and celestial bodies.  

\section*{Acknowledgements}
{Work is supported by T\"{U}B\.{I}TAK under the grant number 114C114.}

\def\cprime{$'$}

\end{document}